\documentstyle[12pt]{article}

\begin{document}
\newcommand{\ol }{\overline}
\newcommand{\ul }{\underline }
\newcommand{\ra }{\rightarrow }
\newcommand{\lra }{\longrightarrow }
\newcommand{\ga }{\gamma }
\newcommand{\st }{\stackrel }
\newcommand{\scr }{\scriptsize }
\title{\Large\textbf{On Nilpotent Multipliers of some Verbal Products of Groups }}
\author{\textbf{Azam Hokmabadi and Behrooz Mashayekhy \footnote{Corresponding author: mashaf@math.um.ac.ir}} \\
Department of Pure Mathematics,\\ Center of Excellence in Analysis on Algebraic Structures,\\
Ferdowsi University of Mashhad,\\
P. O. Box 1159-91775, Mashhad, Iran}
\date{ }

\maketitle
\begin{abstract}
The paper is devoted to finding a homomorphic image for the
$c$-nilpotent multiplier of the verbal product of a family of
groups with respect to a variety ${\mathcal V}$ when ${\mathcal V}
\subseteq {\mathcal N}_{c}$ or ${\mathcal N}_{c}\subseteq
{\mathcal V}$. Also a structure of the $c$-nilpotent multiplier of
a special case of the verbal product, the nilpotent product, of
cyclic groups is given. In fact, we present an explicit formula
for the $c$-nilpotent multiplier of the $n$th nilpotent product of
the group $G= {\bf {Z}}\stackrel{n}{*}...\stackrel{n}{*}{\bf
{Z}}\stackrel{n}{*} {\bf
{Z}}_{r_1}\stackrel{n}{*}...\stackrel{n}{*}{\bf{Z}}_{r_t}$, where
$r_{i+1}$ divides $r_i$ for all $i$, $1 \leq i \leq t-1$, and
$(p,r_1)=1$ for any prime $p$ less than or equal to $n+c$, for all
positive integers $n$, $c$.
\end{abstract}
$2000$ \textit{Mathematics Subject Classification}: 20E34; 20E10; 20F18; 20C25.

\section{Introduction and Motivation}

Let $G=F/R$ be a free presentation of a group $G$. Then the Baer
invariant of $G$ with respect to the variety ${\mathcal N}_c$ of
nilpotent groups of class at most $c \geq 1$, denoted by
${\mathcal N}_cM(G)$, is defined to be
$${\mathcal N}_cM(G)=\frac{R \cap \gamma_{c+1}(F)}{[R,\ _cF]}.$$
${\mathcal N}_cM(G)$ is also called the $c$-nilpotent multiplier
of $G$. Clearly if $c=1$, then ${\mathcal N}_c = {\mathcal A}$ is
the variety of all abelian groups and the Baer invariant of $G$
with respect to this variety is $$ M(G) = \frac{R\cap F'}{[R,
F]},$$ which is the well-known Schur multiplier of $G$.

It is important to find structures for the Schur multiplier and
its generalization, the $c$-nilpotent multiplier, of some famous
products of groups. Determining these Baer invariants of a given
group is known to be very useful for the classification of groups
into isoclinism classes (see [1]).

In 1907, Schur [17], using a representation method, found a
structure for the Schur multiplier of a direct product of two
groups. Also, Wiegold [19] obtained the same result by some
properties of covering groups. In 1979 Moghaddam [13] found a
formula for the $c$-nilpotent multiplier of a direct product of
two groups, where $c+1$ is a prime number or 4. Also, in 1998
Ellis [2] extended the formula for all $c \geq 1$. In 1997 the
second author and Moghaddam [10] presented an explicit formula for
the $c$-nilpotent multiplier of a finite abelian group for any $c
\geq 1$. It is known that the direct product is a special case of
the nilpotent product and we know that regular and verbal products
are generalizations of the nilpotent product.

In 1972, Haebich [6] found a formula for the Schur multiplier of a
regular product of a family of groups. Then the second author [8]
extended the result to find a homomorphic image with a structure
similar to Haebich's type for the $c$-nilpotent multiplier of a
nilpotent product of a family of groups.

In section two, we extend the above result and find a homomorphic
image for the $c$-nilpotent multiplier of a verbal product of a
family of groups with respect to a variety ${\mathcal V}$ when
${\mathcal V} \subseteq {\mathcal N}_{c}$ or ${\mathcal
N}_{c}\subseteq {\mathcal V}$.

A special case of the verbal product of groups whose nilpotent
multiplier has been studied more than others is the nilpotent
product of cyclic groups. In 1992, Gupta and Moghaddam [5]
calculated the $c$-nilpotent multiplier of the nilpotent dihedral
group of class $n$, i.e. $G_n \cong {\bf{Z}}_2 \stackrel{n}{*}
{\bf{Z}}_2$. (Note that in 2001 Ellis [3] remarked that there is a
slip in the statement and gave the correct one.) In 2003,
Moghaddam, the second author and Kayvanfar [14] extended the
previous result and calculated the $c$-nilpotent multiplier of the
$n$th nilpotent product of cyclic groups for $n=$2, 3, 4 under
some conditions. Also, the second author and Parvizi [11, 12]
presented structures for some Baer invariants of a free nilpotent
group that is the nilpotent product of infinite cyclic groups.
Finally the authors and Mohammadzadeh [9] obtained an explicit
formula for the $c$-nilpotent multiplier of the $n$th nilpotent
product of some cyclic groups
$G={\bf{Z}}\stackrel{n}{*}...\stackrel{n}{*}{\bf{Z}}\stackrel{n}{*}
{\bf{Z}}_{r_1}\stackrel{n}{*}...\stackrel{n}{*}{\bf{Z}}_{r_t}$,
where $r_{i+1}$ divides $r_i$ for all $i$, $1 \leq i \leq t-1$,
for $c\geq n$ such that $(p,r_1)=1$ for any prime $p$ less than or
equal to $n$.

In section three, we give an explicit formula for the
$c$-nilpotent multiplier of the above group $G$ when $(p,r_1)=1$
for any prime $p$ less than or equal to $n+c$, for all positive
integers $c,\ n$.

\section{Verbal products}

A group $G$ is said to be a \textit{regular product} of its
subgroups $A_i$, $i \in I$, where $I$ is an ordered set, if the following two conditions hold:\\
i)  $G=\langle A_i | i \in I \rangle$; \\
ii) $A_i \cap \hat{A}_i=1$ for all $i \in I$, where
$\hat{A}_i= \langle A_j | j\in I, j \neq i \rangle$.\\

\hspace{-0.65cm}\textbf{Definition 2.1}. Consider the map
$$ \psi : \prod _{i \in I } ^{\ \ \ *}  A_{i} \rightarrow \prod _{i \in I }
^{\ \ \ \times} A_{i}$$
$$ a_1a_2 \ldots a_n \mapsto (a_1,a_2, \ldots, a_n),$$
which is a natural map from the free product of $\{A_i \}_{i \in
I}$ on to the direct product of $\{A_i \}_{i \in I}$. Clearly its
kernel is the normal closure of
$$ \langle [A_i,A_j] | i,j \in I, i \neq j \rangle $$ in the free
product $A = \prod^*_{i \in I}A_i$. It is denoted by $[A_i^A]$ and
called \textit{the Cartesian subgroup} of the free product (see
[16] for the properties of cartesian subgroups).\\

The following theorem gives a characterization of a regular
product.\\
\hspace{-0.65cm}\textbf{Theorem 2.2} (Golovin 1956 [4]). Suppose
that a group $G$ is generated by a family $\{A_i | i \in I\}$ of
its subgroups, where $I$ is an ordered set. Then $G$ is a regular
product of the $A_i$ if and only if every element of $G$ can be
written uniquely as a product
$$a_1a_2...a_nu, $$ where $1 \neq a_i \in A_{\lambda_i}$,
$\lambda_1<...<\lambda_n$ and $u \in [A_i^G] = \langle
[A_i^G,A_j^G] | i,j \in I, i \neq j \rangle .$\\

\hspace{-0.65cm}\textbf{Definition 2.3}. Let ${\mathcal V}$ be a
variety of groups defined by a set of laws $V$. Then \textit{the
verbal product} of a family of groups $\{A_i\}_{i \in I}$
associated with the variety ${\mathcal V}$ is defined to be
$${\mathcal V}\prod_{i \in I }A_i = \frac{\prod^*_{i \in
I}A_{i}}{V(A) \cap [A_i^A]}.$$

The verbal product is also known as varietal product or simply
${\mathcal V}$-product. If ${\mathcal V}$ is the variety of all
groups, then the corresponding verbal product is the free product;
if ${\mathcal V}={\mathcal A}$ is the variety of all abelian
groups, then the verbal product is the direct product and if
${\mathcal V}={\mathcal N}_c$ is the variety of all nilpotent
groups of class at most $c$, then the verbal product will be the
nilpotent product.

Let $ \{A_{i}|i \in I \}$ be a family of groups and
$$ 1 \rightarrow R _{i} \rightarrow F_{i} \stackrel{\theta_i}{\rightarrow}
 A_{i} \rightarrow 1 $$
be a free presentation for $ A_{i}$. We denote by $\theta$ the
natural homomorphism from the free product $F=\prod^*_{i \in
I}F_{i}$ onto $A=\prod^*_{i \in I}A_{i}$ induced by the
$\theta_i$. Also we assume that the group $G$ is the verbal
product of $\{A_{i}\}_{i \in I }$ associated with the variety
${\mathcal V}$. If $\psi$ is the natural homomorphism from $A$
onto $G$ induced by the identity map on each $A_i$, then we have
the sequence
$$ F= \prod ^{* } _{i \in I } F_{i} \stackrel{\theta}{\rightarrow} A=\prod^{* } _{i \in I
}A_i \stackrel{\psi_v}{\rightarrow} G= {\mathcal V}\prod_{i \in I
}A_i \rightarrow 1.$$

The following notation will be used throughout this section.\\
\hspace{-0.65cm}\textbf{Notation 2.4.}\\
i)  $D_1= \prod_{i \neq j}[R_i, F_j]^F$; \\
ii) $E_c= D_1 \cap \gamma_{c+1}(F)$; \\
iii)$D_c= \prod_{\exists j \ s.t.\ \mu_j \neq i}[R_i, F_{\mu_1}, ...,F_{\mu_c}]^F$; \\
iv) $K_v= V(F)\cap [F_i^F]$;\\
v)  $K_c= \gamma_{c+1}(F)\cap [F_i^F]$.\\

Let $H_v$ be the kernel of $\psi_v$ and $R$ be the kernel of
$\psi_v \circ \theta$. It is clear that $R$ is actually the
inverse image of $H_v$ in $F$ under $\theta$, where $H_v= V(A)\cap
[A_i^A]$ by the definition of the verbal product. Put $H_c=
\gamma_{c+1}(A)\cap [A_i^A]$, then an immediate consequence is the
following lemma.\

\hspace{-0.65cm}\textbf{Lemma 2.5}. With the above notation we have \\
i) $\theta(K_v)=H_v $ and $\theta(K_c)=H_c$;\\
ii) $G=F/R$ and $R=\prod_{i \in I}R_i^F K_v =(\prod_{i \in I}R_i)D_1K_v.$\\

\hspace{-0.65cm}\textbf{Proof. } (i) This follows from the
definition of $\theta$.\\
(ii) It is easy to see that $\ker \theta = \prod_{i \in I}R_i^F$.
On the other hand, since $ \theta(K_v)= \ker \psi_v$, we have $R=
(\ker \theta)K_v=\prod_{i \in I}R_i^F K_v$. Also for all $r \in
R_i$ and $f \in F$, $r^f= r[r,f]$. This implies that $\prod_{i \in
I}R_i^F=\prod_{i \in I}R_i [R_i,F]$. Since $[R_i,F_i]\subseteq
R_i$, $\prod_{i \in I}R_i^F = \prod_{i \in I}R_i D_1$.\\

We now prove some lemmas to compute the $c$-nilpotent multiplier
of $G$.\\
\hspace{-0.65cm}\textbf{Lemma 2.6}. Keeping the above notation we have\\
i)  $[R,\ _cF]= (\prod_{i \in I}[R_i,\ _cF_i])D_c[K_v,\ _cF]$. \\
ii) If $V(F)\subseteq \gamma_{c+1}(F)$, then
$R\cap\gamma_{c+1}(F)=\prod_{i \in I}
(R_i\cap\gamma_{c+1}(F_i))E_c K_v$.\\
iii) If $\gamma_{c+1}(F)\subseteq V(F)$, then
$R\cap\gamma_{c+1}(F)=\prod_{i \in I}
(R_i\cap\gamma_{c+1}(F_i))E_c K_c$.\\

\hspace{-0.65cm}\textbf{Proof. } i)
\begin{eqnarray*}
[R,\ _cF]&=&[\prod_{i \in I}R_i^FK_v,\ _cF]\\
&=&\prod_{i \in I}[R_i,\ _cF]^F[K_v,\ _cF]\\ &=& (\prod_{i \in
I}[R_i,\ _cF_i])D_c[K_v,\ _cF].
\end{eqnarray*}
ii) Let $g\in R\cap\gamma_{c+1}(F)$. Then $g=r_{\lambda_1}...
r_{\lambda_t}dk$ by Lemma 2.5, where $r_{\lambda_i}\in
R_{\lambda_i}$, $d\in D_1$ and $k\in K_v$. Now consider the
natural homomorphism
$$ \varphi: F= \prod ^{\ \ \ * } _{i \in I } F_{i} \rightarrow
\prod ^{\ \ \ \times} _{i \in I } F_{i}.$$ Since $g \in
\gamma_{c+1}(F)$, $\varphi(g)=(r_{\lambda_1}, ...,
r_{\lambda_t})\in \gamma_{c+1}(\prod ^{\times} _{i \in I }
F_{i})=\prod ^{\times} _{i \in I }\gamma_{c+1}(F_{i})$. Therefore
$r_{\lambda_i} \in \gamma_{c+1}(F_{\lambda_i}) \cap
R_{\lambda_i}$ and then $dk \in \gamma_{c+1}(F)\cap [F_i^F]$. Now
since $k \in V(F) \subseteq \gamma_{c+1}(F)$, we have $d \in
\gamma_{c+1}(F)\cap
D_1 =E_c$ and so the result follows.\\
iii) Since $K_c \subseteq K_v$, $\prod_{i \in I}
(R_i\cap\gamma_{c+1}(F_i))E_c K_c \subseteq R\cap\gamma_{c+1}(F)$.
For the reverse inclusion, similar to part (i), $dk \in
\gamma_{c+1}(F)\cap [F_i^F]$. Therefore
$R\cap\gamma_{c+1}(F)\subseteq \prod_{i \in I}
(R_i\cap\gamma_{c+1}(F_i))K_c$. Now the inclusion $E_c \subseteq
K_c$ shows that the equality (iii) holds. $\Box$\\

\hspace{-0.65cm}\textbf{Lemma 2.7}. With the above notation, let
$\varphi_c: F \rightarrow F/E_c$ be the natural homomorphism.
Then $\varphi_c(\prod_{i \in I} (R_i\cap\gamma_{c+1}(F_i))K_v)$
is the direct product of its subgroups  $\varphi_c(K_v)$ and
$\varphi_c(R_i\cap\gamma_{c+1}(F_i))$, $i \in I$.\\

\hspace{-0.65cm}\textbf{Proof.} The Three Subgroups Lemma shows
that
$$[R_i\cap\gamma_{c+1}(F_i), K_v]\subseteq E_c \ \ \ \ for \ all
\ i\in I$$ and $$[R_i\cap\gamma_{c+1}(F_i),
R_j\cap\gamma_{c+1}(F_j)] \subseteq E_c \ \ \ \ for\ all \ i,j \in
I, i\neq j. $$ So we have  $$[\varphi_c(R_i\cap\gamma_{c+1}(F_i)),
\varphi_c(K_v)]=1\ \ \ \ for \ all \ i\in I $$ and
$$[\varphi_c(R_i\cap\gamma_{c+1}(F_i)),
\varphi_c(R_j\cap\gamma_{c+1}(F_j))]=1\ \ \ \ for\ all \ i,j \in
I, i\neq j.$$ Moreover, by Theorem 2.2 we conclude that
$$\varphi_c(R_i\cap\gamma_{c+1}(F_i))\cap (\prod_{i\neq j}\varphi_c(R_j\cap
\gamma_{c+1}(F_j))\varphi_c(K_v))=1.$$ Now the result follows by
the definition  of the direct product. $\Box$\\

\hspace{-0.65cm}\textbf{Lemma 2.8}. With the previous notation, \\
i)  If $V(F) \subseteq \gamma_{c+1}(F)$, then
$\varphi_c(K_v)/\varphi_c([K_v,\ _cF]) \cong H_v/[H_v,\ _cA]$.\\
ii) If $\gamma_{c+1}(F)\subseteq V(F)$, then
$\varphi_c(K_c)/\varphi_c([K_c,\ _cF] \cong H_c/[H_v,\ _cA]$.\\

\hspace{-0.65cm}\textbf{Proof.} i) If $V(F)\subseteq
\gamma_{c+1}(F)$, then
$$\frac{\varphi_c(K_v)}{\varphi_c([K_v,\ _cF])} \cong
\frac{K_vE_c}{[K_v,\ _cF]E_c}\cong \frac{K_v}{K_v \cap [K_v,\
_cF]E_c}.$$ On the other hand
$$\frac{\theta(K_v)}{\theta([K_v,\ _cF])} \cong
\frac{K_v \ker\theta}{[K_v,\ _cF]\ker\theta}\cong \frac{K_v}{K_v
\cap [K_v,\ _cF]\ker\theta}\cong \frac{K_v}{K_v \cap [K_v,\
_cF]D_1\prod_{i \in I}R_i}.$$ Now Theorem 2.2 and definition of
$E_c$ imply that $$\frac{\theta(K_v)}{\theta([K_v,\ _cF])} \cong
\frac{K_v}{K_v \cap [K_v,\ _cF]E_c}.$$ Therefore by Lemma 2.5, we
conclude that
$$\frac{\varphi_c(K_v)}{\varphi_c([K_v,\ _cF])} \cong \frac{\theta
(K_v)}{\theta([K_v,\ _cF])} \cong \frac{H_v}{[H_v,\ _cA]}.$$
ii) The proof is similar to (i).\\

Now we are ready to state and prove the main result of this
section. \\
\hspace{-0.65cm}\textbf{Theorem 2.9}. With the above notation, \\
i) If ${\mathcal N}_{c}\subseteq {\mathcal V}$, then
$\prod^\times _{i \in I}{\mathcal N}_{c}M(A_i) \times H_v/[H_v,\
_cA]$ is a homomorphic
image of ${\mathcal N}_{c}M({\mathcal V}\prod_{i \in I}A_i)$, and if
${\mathcal V}\prod_{i \in I}A_i$ is finite, then the above structure is isomorphic
to a subgroup of ${\mathcal N}_{c}M({\mathcal V}\prod_{i \in I}A_i)$.\\
ii) If ${\mathcal V}\subseteq {\mathcal N}_{c}$, then
$\prod^\times _{i \in I}{\mathcal N}_{c}M(A_i) \times H_c/[H_v,\
_cA]$ is a homomorphic image of ${\mathcal N}_{c}M({\mathcal
V}\prod_{i \in I}A_i)$, and if ${\mathcal V}\prod_{i \in I}A_i$
is finite, then the above structure is isomorphic
to a subgroup of ${\mathcal N}_{c}M({\mathcal V}\prod_{i \in I}A_i)$.\\

\hspace{-0.65cm}\textbf{Proof.}  i) By Lemma 2.6 (i),(ii)
$${\mathcal N}_{c}M({\mathcal V}\prod_{i \in I}A_i)\cong \frac{R
\cap \gamma_{c+1}(F)}{[R,\ _cF]}\cong\frac{\prod_{i \in I}
(R_i\cap\gamma_{c+1}(F_i))E_c K_v}{\prod_{i \in I}[R_i,\
_cF_i]D_c[K_v,\ _cF]}.$$ Therefore there is a natural epimorphism
from ${\mathcal N}_{c}M({\mathcal V}\prod_{i \in I}A_i)$ to $$
\frac{\prod_{i \in I} (R_i\cap\gamma_{c+1}(F_i))E_c K_v}{\prod_{i
\in I}[R_i,\ _cF_i]E_c[K_v,\ _cF]} \cong \frac{\varphi_c(\prod_{i
\in I} (R_i\cap\gamma_{c+1}(F_i)) K_v)}{\varphi_c(\prod_{i \in
I}[R_i,\ _cF_i][K_v,\ _cF])}.$$ Lemma 2.7 and the fact that
$\varphi_c([K_v,\ _cF])\subseteq \varphi_c(K_v)$ and
$\varphi_c([R_i,\ _cF_i])\subseteq \varphi_c(
R_i\cap\gamma_{c+1}(F_i))$ imply that
$$\frac{\varphi_c(\prod_{i\in I} (R_i\cap\gamma_{c+1}(F_i))
K_v)}{\varphi_c(\prod_{i \in I}[R_i,\ _cF_i][K_v,\ _cF])} \cong
\prod^\times_{i \in I}\frac{\varphi_c(
R_i\cap\gamma_{c+1}(F_i))}{\varphi_c([R_i,\ _cF_i])} \times
\frac{\varphi_c(K_v)}{\varphi_c([K_v,\ _cF])}.$$ It is
straightforward to see that
$$\frac{\varphi_c( R_i\cap\gamma_{c+1}(F_i))}{\varphi_c([R_i,\ _cF_i])} \cong
\frac{R_i\cap\gamma_{c+1}(F_i)}{[R_i,\ _cF_i]}$$ by Theorem 2.2.
Therefore, the result holds by Lemma 2.8 (i).\\
ii) By an argument similar to (i), we obtain the result. $\Box$\\

We need the following lemma whose proof is straightforward.\\
\hspace{-0.65cm}\textbf{Lemma 2.10}. Let $ \{A_{i}|i \in I \}$ be
a family of groups. Put $A=\prod^*_{i \in I}A_i$. Then for all
integers $m \geq 2$,
$$\gamma_m(A) = \prod_{i \in I}\gamma_m(A_i)(\gamma_m(A)\cap
[A_i^A]).$$ In particular if the $A_i$ are cyclic, then
$\gamma_m(A) = \gamma_m(A)\cap [A_i^A]$. \\

The following corollary is an interesting consequence of Theorem
2.9 for cyclic groups.\\
\hspace{-0.65cm}\textbf{Corollary 2.11}. Let $ \{A_{i}|i \in I \}$
be a family of cyclic groups. Then \\
i) If ${\mathcal N}_{c}\subseteq {\mathcal V}$, then ${\mathcal
N}_{c}M({\mathcal V}\prod_{i \in I}A_i) \cong H_v/[H_v,\ _cA]$.
Moreover if ${\mathcal V}\subseteq {\mathcal N}_{2c}$, then
$V(\prod _{i \in I } ^{\stackrel{2c}{*}}A_{i})$ is a homomorphic
image of ${\mathcal N}_{c}M({\mathcal V}\prod_{i \in I}A_i)$.\\
ii) If ${\mathcal V}\subseteq {\mathcal N}_{c}$, then ${\mathcal
N}_{c}M({\mathcal V}\prod_{i \in I}A_i) \cong H_c/[H_v,\ _cA]$.
Moreover if ${\mathcal N}_{m}\subseteq {\mathcal V}$, then
$\gamma_{c+1}(\prod _{i \in I } ^{\stackrel{m+c}{*}}A_{i})$ is a
homomorphic image of ${\mathcal N}_{c}M({\mathcal V}\prod_{i \in I}A_i)$.\\

\hspace{-0.65cm}\textbf{Proof.} i) Since the $A_i$ are cyclic
groups and the $R_i$ have no commutators, it is concluded that
$D_c=E_c$. So the epimorphism in the proof of Theorem 2.9, is
actually an isomorphism. Also ${\mathcal N}_{c}M(A_i)=1$,
therefore ${\mathcal N}_{c}M({\mathcal V}\prod_{i \in I}A_i)
\cong H_v/[H_v,\ _cA]$. Now suppose ${\mathcal N}_c \subseteq
{\mathcal V}\subseteq {\mathcal N}_{2c}$. The inclusion
$V(A)\subseteq \gamma_{c+1}(A)$ and Lemma 2.10 imply that $V(A)
\subseteq [A_i^{A}]$ and thus $H_v=V(A)\cap [A_i^{A}]=V(A)$. So
we have ${\mathcal N}_{c}M({\mathcal V}\prod_{i \in I}A_i)=
V(A)/[V(A),\ _cA]$ and hence $V(A)/ \gamma_{2c+1}(A)$ is a
homomorphic image of ${\mathcal N}_{c}M({\mathcal V}\prod_{i \in
I}A_i)$. On the other hand since ${\mathcal V}\subseteq {\mathcal
N}_{2c}$, we have $V(A)/ \gamma_{2c+1}(A)=V(A/ \gamma_{2c+1}(A))=
V(\prod _{i \in I } ^{\stackrel{2c}{*}}A_{i})$. This completes the proof.\\
ii) An argument similar to (i), shows that ${\mathcal
N}_{c}M({\mathcal V}\prod_{i \in I}A_i) \cong H_c/[H_v,\ _cA]$.
Now since ${\mathcal N}_{m}\subseteq {\mathcal V} \subseteq
{\mathcal N}_c$, $\gamma_{c+1}(A)/\gamma_{m+c+1}(A)$ is a
homomorphic image of ${\mathcal N}_{c}M({\mathcal V}\prod_{i \in
I}A_i)$ and also
$$\frac{\gamma_{c+1}(A)}{\gamma_{m+c+1}(A)}=\gamma_{c+1}(\frac{A}{\gamma_{m+c+1}(A)})
=\gamma_{c+1}(\prod_{i \in I } ^{\stackrel{m+c}{*}}A_{i}).$$ Hence
the result follows. $\Box$\\

\hspace{-0.65cm}\textbf{Remark 2.12}.  Let $ \{A_{i}|i \in I \}$
be a family of groups. \\
i) If ${\mathcal V}$ is the variety of trivial groups, then
Theorem 2.9 implies that $\prod^\times _{i \in I}{\mathcal
N}_{c}M(A_i)$ is a homomorphic image of ${\mathcal
N}_{c}M(\prod^*_{i \in I}A_i)$. In particular $M(\prod^*_{i \in
I}A_i)= \prod^\times _{i \in I}M(A_i)$
which is a result of Miller [15].\\
ii) If ${\mathcal V}$ is the variety of nilpotent groups of class
at most $n$, ${\mathcal N}_n$, then main results of the second
author [8] are obtained by Theorem 2.9 and corollary 2.11.

\section{Nilpotent Products of Cyclic Groups}

In this section we use a result of the previous section and find a
structure for the $c$-nilpotent multiplier of the group
$G={\bf{Z}}\stackrel{n}{*}...\stackrel{n}{*}{\bf{Z}}\stackrel{n}{*}
{\bf{Z}}_{r_1}\stackrel{n}{*}...\stackrel{n}{*}{\bf{Z}}_{r_t}$,
where $r_{i+1}$ divides $r_i$ for all  $i$, $1 \leq i \leq t-1$,
 such that $(p,r_1)=1$ for any prime $p$ less than or equal to
$n+c$. The proof relies on basic commutators [7] and related
results. We recall that the number of basic commutators of weight
$c$ on $n$ generators, denoted by $\chi_{c}(n)$, is determined by
Witt formula [7]. Also, M. Hall proved that if $F$ is the free
group on free generators $x_{1}, x_{2}, ...,x_{r}$ and $c_1, ...,
c_t$ are basic commutators of weight $1, 2, ..., n$, on $x_1,
\dots ,x_r$, then an arbitrary element $f$ of $F$ has a unique
representation,$$ f = c_1^{\beta_1}c_2^{\beta_2} ...
c_t^{\beta_t} \ \ mod \gamma_{n+1}(F).$$ In particular the basic
commutators of weight $n$ provide a basis for the free abelian
group
$\gamma_n(F)/\gamma_{n+1}(F)$ (see [7]).\\

The following theorem represents the elements of some nilpotent
products of cyclic groups in terms of basic commutators.

\hspace{-0.65cm}\textbf{Theorem 3.1} ([18]). Let $A_1,...,A_t$ be
cyclic groups of order $\alpha_1, ...,\alpha_t$ respectively,
where if $A_i$ is infinite cyclic, then $\alpha_i=0$. Let $a_i$
generate $A_i$ and let $G= A_1 \stackrel{n}{*}  ...
\stackrel{n}{*} A_{t}$, where $n$ is greater than or equal to 2.
Suppose that all the primes appearing in the factorizations of
the $\alpha_i$ are greater than or equal to $n$ and $u_1,u_2,...,$
are basic commutators of weight less than $n$, on the letters
$a_1, \dots, a_t$. Put $N_i=\alpha_{i_j}$ if $u_i=a_{i_j}$ of
weight 1, and
$$ N_i = gcd(\alpha_{i_1}, ...,\alpha_{i_k}) $$ if $a_{i_j}$, $1
\leq j \leq k$, appears in $u_i$. Then every element $g$ of $G$
can be uniquely expressed as $$g=\prod u_i^{m_i},$$ where the
$m_i$ are integers modulo $N_i$ (by $gcd$ we mean the greatest
common divisor).\\

The following theorem is an interesting consequence of Corollary
2.11.

\hspace{-0.65cm}\textbf{Theorem 3.2}. Let $ \{A_{i}|i \in I \}$
be a family of cyclic groups. Then \\
i) if $n \geq c$, then ${\mathcal
N}_{c}M(\prod^{\stackrel{n}{*}}_{i \in I}A_i)\cong
\gamma_{n+1}(\prod^{\stackrel{n+c}{*}}_{i \in I}A_i)$;\\
ii) if $c \geq n$, then ${\mathcal
N}_{c}M(\prod^{\stackrel{n}{*}}_{i \in
I}A_i)\cong \gamma_{c+1}(\prod^{\stackrel{n+c}{*}}_{i \in I}A_i)$.\\

\hspace{-0.65cm}\textbf{Proof.} i) Put ${\mathcal V}= {\mathcal
N}_{n}$ in Corollary 2.11 and deduce that
$${\mathcal N}_{c}M(\prod^{\stackrel{n}{*}}_{i \in I}A_i) \cong
H_n/[H_n,\ _cA].$$ On the other hand by Lemma 2.10, $H_n=
\gamma_{n+1}(A)\cap [A_i^A]= \gamma_{n+1}(A)$. Therefore
$${\mathcal N}_{c}M(\prod^{\stackrel{n}{*}}_{i \in
I}A_i)  \cong  \frac{\gamma_{n+1}(A)}{[\gamma_{n+1}(A),\ _cA]}=
\gamma_{n+1}(\frac{A}{\gamma_{n+c+1}(A)}) =
\gamma_{n+1}(\prod^{\stackrel{n+c}{*}}_{i \in I}A_i).$$
ii) The result follows as for (i). $\Box$\\

Now, we are in a position to state and prove the main result of
this section .\\
\hspace{-0.65cm}\textbf{Theorem 3.3}. Let $
G=A_1\stackrel{n}{*}...\stackrel{n}{*}A_{m+t}$ be the $n$th
nilpotent product of cyclic groups such that $A_i\cong {\bf Z}$
for $1\leq i \leq m$ and $A_{m+j}\cong {\bf Z}_{r_j}$ and $r_{j+1}
\mid r_j$ for all $ 1 \leq j \leq t-1$. If
$(p,r_1)=1$ for any prime $p$ less than or equal to $n+c$, then \\
i) if $n \geq c$, then ${\mathcal N}_{c}M(G) \cong
{\bf{Z}}^{(g_0)} \oplus {\bf{Z}}_{r_1}^{(g_1 - g_0)}\oplus
...\oplus {\bf{Z}}_{r_t}^{(g_t - g_{t-1})}$;\\
ii) if $c \geq n$, then ${\mathcal N}_{c}M(G) \cong
{\bf{Z}}^{(f_0)} \oplus {\bf{Z}}_{r_1}^{(f_1 - f_0)}\oplus
...\oplus {\bf{Z}}_{r_t}^{(f_t - f_{t-1})}$,\\
where $f_k=\sum ^{n}_{i=1} \chi _{c+i}(m+k)$ and $g_k=\sum
^{c}_{i=1} \chi _{n+i}(m+k)$ for $0 \leq k \leq t$ and
${\bf{Z}}_{r}^{(d)}$ denotes the direct sum of $d$ copies of the
cyclic group ${\bf{Z}}_{r}$.\\

\hspace{-0.65cm}\textbf{Proof.} i) If $n \geq c$, then by Theorem
3.2, it is enough to find the structure of
$\gamma_{n+1}(\prod^{\stackrel{n+c}{*}}_{i \in I}A_i)$. Suppose
that $a_i$ generates $A_i$ and $F$ is the free group generated by
$a_1, ...,a_{m+t}$. Let $B$ be the set of all basic commutators of
weight $1,2,..., c+n$ on the letters $a_1, ...,a_{m+t}$. Now
define
$$D=\{ u^{N_i}| \ u \in B \ and \ N_i=gcd(\alpha_{i_1},
...,\alpha_{i_k}) \ if \ a_{i_j} \ appears \ in \ u \ for \ 1
\leq j \leq k \}.$$ Then Theorem 3.1 implies that
$\prod^{\stackrel{n+c}{*}}_{i \in I}A_i= F/\langle D
\rangle\gamma_{c+n+1}(F)$ and so
\begin{eqnarray*}
\gamma_{n+1}(\prod^{\stackrel{n+c}{*}}_{i \in I}A_i) &=&
\gamma_{n+1} (\frac{F}{\langle D \rangle \gamma_{c+n+1}(F)})\\ &=&
\frac{\gamma_{n+1}(F)}{\langle D \rangle\gamma_{c+n+1}(F) \cap
\gamma_{n+1}(F)}\\ &\cong&
\frac{\gamma_{n+1}(F)/\gamma_{c+n+1}(F)}{(\langle D \rangle \cap
\gamma_{n+1}(F))\gamma_{c+n+1}(F)/ \gamma_{c+n+1}(F)}.
\end{eqnarray*}

It can be deduced from Hall Theorem that
$\gamma_{n+1}(F)/\gamma_{c+n+1}(F)$ is a free abelian group with
a basis $\bar{B_1}=\{ u \gamma_{c+n+1}(F) | \ u \in B_1 \}$,
where $B_1$ is the set of all basic commutators of weight
$n+1,..., c+n$ on $a_1, ...,a_{m+t}$. Also, the uniqueness of the
presentation of elements implies that the abelian group $(\langle
D \rangle \cap \gamma_{n+1}(F))\gamma_{c+n+1}(F)/
\gamma_{c+n+1}(F)$ is free with a basis
$$\bar{E}=\{ u \gamma_{c+n+1}(F) | u \in D \cap
\bar{B_1}= \bigcup_{j=1}^tD_j \} \ ,$$ where $D_j$ is the set of
all $u^{r_j}$, such that $u$ is a basic commutator of weight
$n+1,..., c+n$ on $a_1, ...,a_{m+j}$ such that $a_{m+j}$ appears
in $u$. Also we have $$|D_j|= \sum_{i=1}^c \chi _{n+i}(m+j)- \chi
_{n+i}(m+j-1)= g_j -g_{j-1}.$$ This
completes the proof. \\
ii) The proof is similar to (i). $\Box$ \\

Note that the authors with F. Mohammadzadeh [9] by a different
method presented a similar structure for ${\mathcal N}_{c}M(G)$,
for $c\geq n$ with a weaker condition $(p,r_1)=1$ for any prime
$p$ less than or equal to $n$.\\

\hspace{-0.65cm}\textbf{Remark 3.4}. The condition $r_{j+1} \mid
r_j$, in the above theorem, simplifies the structure of the
$c$-nilpotent multiplier of $G$ and gives a clear formula. One can
use the above method and find the structure of ${\mathcal
N}_{c}M(G)$ without the condition $r_{j+1} \mid r_j$, but with a
more complex formula. For example, for a simple case if $
G={\bf{Z}}_r \stackrel{n}{*} {\bf{Z}}_s$ where $(p,r)=(p,s)=1$ for
any prime $p$ less than or equal to $n+c$ and $(r,s)=d$, then \\
i) if $n \geq c$, then ${\mathcal N}_{c}M(G) \cong
{\bf{Z}}_{d}^{(\sum ^{c}_{i=1} \chi _{n+i}(2))}$;\\
ii) if $c \geq n$, then ${\mathcal N}_{c}M(G) \cong
{\bf{Z}}_{d}^{(\sum ^{n}_{i=1} \chi _{c+i}(2))}$.\\

\end{document}